\documentclass{amsart}

\usepackage{amsmath, amsthm, amssymb}
\usepackage{hyperref}
\usepackage{mathbbol}
\usepackage{graphicx}
\usepackage{color}

\newtheorem{thm}{Theorem}

\newtheorem{lemma}[thm]{Lemma}
\theoremstyle{definition}
\newtheorem{defi}[thm]{Definition}
\newtheorem{cor}[thm]{Corollary}

\theoremstyle{remark}

\newtheorem{rem}[thm]{Remark}

\newcommand{\alg}[1]{\mathfrak{{#1}}}

\newcommand{\cn}[1]{\mathbb{C}^{{#1}}}
\newcommand{\rn}[1]{\mathbb{R}^{{#1}}}
\newcommand{\co}[2]{\left[{#1},{#2}\right]} 

\newcommand{\eref}[1]{(\ref{#1})} 
\newcommand{\vecf}[1]{{ \mathfrak{X} {#1} }} 

\newcommand{\edge}{{\rightarrow }} 
\newcommand{\morphU}{{\mathcal{U} }} 
\newcommand{\widebar}[1]{{\overline{#1}}}

\newcommand{\Hom}{\mathop{Hom}}

\newcommand{\C}{{\mathbb{C}}}
\newcommand{\R}{{\mathbb{R}}}

\newcommand{\gf}{\C}
\newcommand{\tW}{ {[[u]]\otimes_{\gf[u]} W} }  

\begin{document}
\title{Formality of Cyclic Chains}

\author{Thomas Willwacher}

\address{Department of Mathematics, ETH Zurich}
\email{thomas.willwacher@math.ethz.ch}
\thanks{The author was partially supported by the Swiss National Science Foundation (grant 200020-105450)}
\date{}

\begin{abstract}
We prove a conjecture raised by Tsygan \cite{tsygan}, namely the existence of an $L_\infty$-quasiisomorphism of $L_\infty$-modules between the cyclic chain complex of smooth functions on a manifold and the differential forms on that manifold. Concretely, we prove that the obvious $u$-linear extension of Shoikhet's morphism of Hochschild chains solves Tsygan's conjecture.
\end{abstract}
\maketitle

\section{Introduction and notations}
Let $M$ be a smooth manifold and $\vecf{M}$ the Lie algebra of vector fields on $M$. Let $T_{poly}^\bullet=\bigwedge^{\bullet+1} \vecf{M}$ be the algebra of polyvector fields on $M$. It is naturally endowed with a Lie bracket $\co{\cdot}{\cdot}_{SN}$, the Schouten-Nijenhuis bracket. Denote by $A=C^\infty(M)$ the commutative algebra of smooth functions on $M$. Let $D_{poly}$ be the subcomplex of the Hochschild complex $C^\bullet(A,A)$ given by polydifferential operators. The $n$-cochains in this complex are spanned by maps of the form
\[
A^{\otimes n} \ni a_1\otimes\dots \otimes a_n \mapsto \prod_{k=1}^n(D_ka_k) \in A
\]
where the $D_k$ are differential operators. The Hochschild differential $d_H$ and the Gerstenhaber bracket $\co{\cdot}{\cdot}_{G}$ naturally restrict to the subcomplex $D_{poly}$ and endow it with the structure of a differential graded Lie algebra.

In his famous paper \cite{kontsevich} Kontsevich proved in 1997 the \emph{Formality Theorem} (on cochains), i.e., the existence of an $L_\infty$-quasiisomorphism of differential graded Lie algebras
\[
T_{poly} \rightarrow D_{poly}.
\]

The Taylor coefficients of this morphism were explicitly given in terms of graphs. Kontsevich's techniques for dealing with graphs and generating proofs based on Stokes' Theorem are very relevant for most papers on the subject, and the present one is no exception. However, we will not review his construction here, but refer the reader to the original work \cite{kontsevich}.

Next consider the Hochschild chain complex $C_\bullet(A,A)$ of $A$ with values in $A$. It forms a (dgla) module over the cochain complex $C^\bullet(A,A)$, where the action is given by 
\begin{multline}
\label{equ:hochaction}
D \cdot (a_0\otimes \cdots \otimes a_n) =
\sum_{j=n-d+1}^{n}(-1)^{n(j+1)}  D(a_{j+1}, \dots, a_0, \dots) \otimes a_{d+j-n}\otimes \cdots \otimes a_j
+\\
+\sum_{i=0}^{n-d}(-1)^{(d-1)(i+1)}a_0 \otimes \cdots \otimes a_i \otimes D(a_{i+1},\dots, a_{i+d}) \otimes \cdots \otimes a_n
\end{multline}
for $D \in C^d(A,A)$ and $a_0, \dots, a_n \in A$.

Through Kontsevich's morphism the chains $C_\bullet(A,A)$ also carry an $L_\infty$-module structure over the dgla $T_{poly}$.

Furthermore, there is another natural module over $T_{poly}$ that can be constructed without additional data, namely the differential forms $\Omega^\bullet(M)$, with the action given by Lie derivatives
\[
\gamma \cdot \omega = (d \iota_\gamma - (-1)^{p}\iota_\gamma d) \omega
\]
where $\gamma \in \bigwedge^p \vecf{M}\subset T_{poly}$ and $\omega\in \Omega^\bullet(M)$.

A natural extension of the formality Theorem is then the following statement, which was conjectured by Tsygan \cite{tsygan} in 1999.
\begin{thm}[Formality Theorem on Chains]
\label{thm:chthm}
There exists an $L_\infty$-quasiisomorphism of $L_\infty$ modules over $T_{poly}$
\[
\morphU:(C_\bullet(A,A), b) \rightarrow  (\Omega^\bullet (M), 0)
\]
\end{thm}
Here the notation means that the complex $\Omega^\bullet (M)$ is endowed with 0 differential. 
The Theorem has been proven independently by Shoikhet \cite{shoikhet} and Dolgushev \cite{dolgushev} and by Tamarkin and Tsygan \cite{tamarkin}. More precisely, Shoikhet found an explicit quasiisomorphism $\morphU^{sh}$ in the cases $M=\R^n$, or $M$ a formal completion of $\R^n$ at the origin. Dolgushev globalized this construction using Fedosov resolutions. 
The explicit construction of $\morphU^{sh}$ given by Shoikhet will be reviewed in section \ref{sec:shoi}.

Tsygan also conjectured the analog of the above theorem on cyclic instead of Hochschild chains. This is the conjecture that will be proven in the present paper. There are several variants of the cyclic chain complex, all of which have the form
\begin{equation}
\label{equ:cyccdef}
CC^W_p(A) = (C_\bullet(A,A) \tW)_p
\end{equation}
where $W$ is a module over the graded algebra $\cn{}[u]$, with $u$ being a formal variable of degree $-2$.\footnote{This notation is due to Getzler.}
The differential on the above complexes is given by $b+uB$, where $b$ is the Hochschild boundary operator and $B$ is defined by
\begin{equation}
\label{equ:Bdef}
B (a_0\otimes \dots \otimes a_n) = \sum_{j=0}^n (-1)^{nj} 1 \otimes a_j \otimes \dots \otimes a_n\otimes a_0\otimes \dots \otimes a_{j-1}.
\end{equation}
where $a_{-1}:=a_n$ to simplify notation. The homology $HC^W_\bullet (A)$ of the cyclic chain complex is related to the de Rham cohomology of $M$ via the following theorem, which can be found in \cite{block} (Theorem 3.3 for $G=\{pt\}$).
\begin{thm}
\label{thm:hkr}
Let $W$ be a $\cn{}[u]$-module of finite projective dimension over $\C[u]$, then
\[
HC^W_\bullet(A) \cong H^\bullet(\Omega(M) \tW , u d).
\]
\end{thm}

We will prove the following 
\begin{thm}
\label{thm:main}
Shoikhet's $L_\infty$-morphism $\morphU^{sh}$ satisfies
\[
\morphU^{sh} \circ B = d \circ \morphU^{sh}
\]
\end{thm}

As a corollary, one obtains the formality theorem on cyclic chains.
\begin{cor}
\label{cor:tsy}
For $W$ a $\C[u]$-module of finite projective dimension over $\C[u]$, there is an $L_\infty$-quasiisomorphism of $L_\infty$-modules over $T_{poly}$
\[
\morphU: (CC^W_\bullet(A,A),b+uB) \rightarrow (\Omega(M) \tW , u d)
\]
\end{cor}
\begin{proof}
For the proof, one needs to consider Fedosov resolutions of the above two complexes. Introducing these and the required notations would be very lengthy. To avoid this, we take the liberty to copy the notations of Dolgushev, as used in section 5 of \cite{dolgushev}, until the end of this proof. For definitions and explanations, we refer to Dolgushev's diligent treatment. Concretely, there is the following sequence of quasiisomorphisms of $L_\infty$-modules over $T_{poly}$:
\[
C^{poly}(M) \xrightarrow{\rho}  (\Omega(M,\mathcal{C}^{poly}),D+\mathfrak{b}) 
\xrightarrow{\mathfrak{K}}
(\Omega(M,\mathcal{E}),D) \xleftarrow{\tau} \mathcal{A}^\bullet(M).
\]
From left to right, the objects are the Hochschild chain complex of $C^\infty(M)$, its Fedosov resolution, the Fedosov resolution of the de Rham complex and the de Rham complex itself. The middle quasiisomorphism (i.e., $\mathfrak{K}$) is defined using Shoikhet's morphism $\morphU^{sh}$ fiberwise.

All the above four complexes are, in fact, \emph{mixed complexes}, in the sense that they carry another differential of degree $+1$, anticommuting with their boundary operators. This differential is (from left to right) Connes' $B$ as in \eref{equ:Bdef}, the same operator applied fiberwise $B_f$, the fiberwise de Rham differential $d_f$ and finally the de Rham differential $d$. We claim that all morphisms in the above sequence are morphisms of mixed complexes, i.e., commute with the application of the additional differentials. For the middle morphism $\mathfrak{K}$, this follows from Theorem \ref{thm:main} above. For the left- and rightmost morphisms, note that the fiberwise $B_f$ and $d_f$ map $D$-constant sections to $D$-constant sections. Hence it suffices to observe that for $s\in C^{poly}(M)$, $\alpha \in \mathcal{A}^\bullet(M)$, the parts of degree 0 in the formal variable (usually called ``$y$'') of $B_f\rho(s)$ and $d_f\tau(\alpha)$ agree with $Bs$ and $d\alpha$ respectively. 

By $u$-linear extension and Remark \ref{rem:whattoproof} in the appendix, we then obtain the following sequence of morphisms of $L_\infty$-modules over $T_{poly}$:
\begin{multline*}
(C^{poly}(M)\tW,\mathfrak{b}+ uB) \to  (\Omega(M,\mathcal{C}^{poly})\tW,D+\mathfrak{b}+uB_f) \to \\
\to 
(\Omega(M,\mathcal{E})\tW,D+ud_f) \leftarrow (\mathcal{A}^\bullet(M)\tW, ud).
\end{multline*}

It remains to be shown that all these morphisms are quasiisomorphisms. For this, one can forget about the higher degree Taylor components of the $L_\infty$-module-morphisms and consider the above sequence as a sequence of morphisms of complexes. But we know that the (0-th Taylor components of the) original morphisms $\rho$, $\mathfrak{K}$ and $\tau$ were morphisms of mixed complexes inducing isomorphisms on homology (wrt. the degree -1 differential). Hence Proposition 2.4 of \cite{gj} finishes the proof of the Theorem.

\end{proof}

\subsection{Structure of the Paper}
The precise definitions of structures, brackets, differentials and gradings that were omitted in the introduction can be found in the appendix. The author wants to avoid having the reader browse through pages of definitions she or he already knows. So in the next section we directly start by reviewing the construction of Shoikhet's formality morphism, adding several remarks that will simplify the proof of Theorem \ref{thm:main}. The proof can then be found in section \ref{sec:theproof}.

\subsection{Acknowledgements}
The author is grateful to his advisor Prof. Giovanni Felder for introducing him to the problem and many helpful discussion and corrections to this manuscript.

\section{Shoikhet's Formality Theorem on Chains}
\label{sec:shoi}
In this section we recall the construction of Shoikhet's morphism $\morphU^{sh}$ for the case $M=\rn{d}$ and outline his proof of Theorem \ref{thm:chthm}. As usual in deformation quantization, the morphism can be expressed as a sum of graphs. Denote by $\morphU^{sh}_m$ the $m$-th Taylor component of $\morphU^{sh}$. For $\xi$ a constant polyvector field, we will set
\[
\morphU^{sh}_m(\gamma_1\wedge\dots \wedge\gamma_m;a_0\otimes\dots\otimes a_n)[\xi]=\sum_{\Gamma\in G(m,n)} w_\Gamma D_\Gamma(\xi, \gamma_1\wedge\dots \wedge\gamma_m ; a_0,..,a_n).
\]
Here the sum is over all Kontsevich graphs with $m+1$ type I and $n+1$ type II vertices. The polydifferential operator $D_\Gamma$ is the same as in the Kontsevich case, but with the polyvector field $\xi$ put exclusively at the first vertex of the graph. However, the weight $w_\Gamma\in \cn{}$ is defined differently, a formula will be given below. The square brackets shall denote evaluation of a differential form on a polyvector field.

To be precise, we will use here the following definition of the graphs occuring in the sum.

\begin{defi}
The set $G(n,m)$, $n,m\in \mathbb{N}_0$ consists of directed graphs $\Gamma$ such that
\begin{itemize}
\item The vertex set of $\Gamma$ is 
\[
V(\Gamma) = \{0,1,..,n\}\cup \{\bar{0},..,\bar{m}\}
\]
where the vertex $0$ will be called the \emph{central} vertex, the vertices $\{0,1,..,n\}$ the \emph{type I} vertices and the $\{\bar{0},..,\bar{m}\}$ the \emph{type II} vertices.
\item Every edge $e=(v_i\edge v_j)\in E(\Gamma)$ starts at a type I vertex and does not end at the central vertex. I.e., $v_i$ is type I and $v_j$ is not the vertex $0$. We will call the edges $(0\edge v_k)$ that start at the central vertex \emph{central edges} and denote the set of these edges by $E_c(\Gamma)$.
\item There are no tadpoles, i.e., no edges of the form $(v\edge v)$.
\item For each type I vertex $v$, there is an ordering given on
\[
Star(v) = \{(v\edge w) \mid (v\edge w)\in E(\Gamma),\; w\in E(\Gamma) \}.
\] 
\end{itemize}
\end{defi}

Let us next define the weight $w_\Gamma$ of $\Gamma\in G(n,m)$. As in the Kontsevich case, it is an integral of a certain differential form over a compact manifold with corners, the \emph{configuration space} $C_\Gamma$.
\begin{equation}
\label{equ:wgammadef}
w_\Gamma = \prod_{v\in V(\Gamma)} \frac{1}{(\#Star(v))!} \int_{C_\Gamma} \omega_\Gamma
\end{equation}

\begin{defi}
\label{def:configspace}
The {enlarged configuration space} $\tilde{C}_\Gamma$ is the Fulton-MacPherson-like\footnote{We mean the compactification constructed similarly to \cite{kontsevich}, section 5. It will not be of any importance.} compactification of the space of embeddings 
\[
(z_0,\dots,z_n,z_{\bar{0}},\dots z_{\bar{m}}): V(\Gamma) \rightarrow D
\]
of the vertex set $V(\Gamma)$ of $\Gamma$ into the closed unit disk $D=\{z \in \cn{}; |z|\leq 1\}$ such that
\begin{enumerate}
\item The central vertex is mapped to the origin, i.e., $z_0=0$.
\item All type I vertices are mapped to the interior of $D$, i.e. $z_j\in D^\circ$ for $j=1,..,n$.
\item All type II vertices are mapped to the boundary of $D$, i.e. $z_{\bar{j}}\in \partial D$ for $j=0,..,m$.
\item The type II vertices occur in counterclockwise increasing order on the circle, i.e., $0<\arg \frac{z_{\bar{1}}}{z_{\bar{0}}}< \dots < \arg \frac{z_{\bar{n}}}{z_{\bar{0}}}<2\pi$.
\end{enumerate}
The \emph{configuration space} $C_\Gamma$ is the codimension 1 subspace of $\tilde{C}_\Gamma$ on which $z_{\bar{0}}=1$, i.e., $C_\Gamma = \{z_{\bar{0}}=1\}\subset \tilde{C}_\Gamma$
\end{defi}

\begin{figure}
\includegraphics[trim = 0mm 60mm 0mm 20mm, clip, width=.5\textwidth]{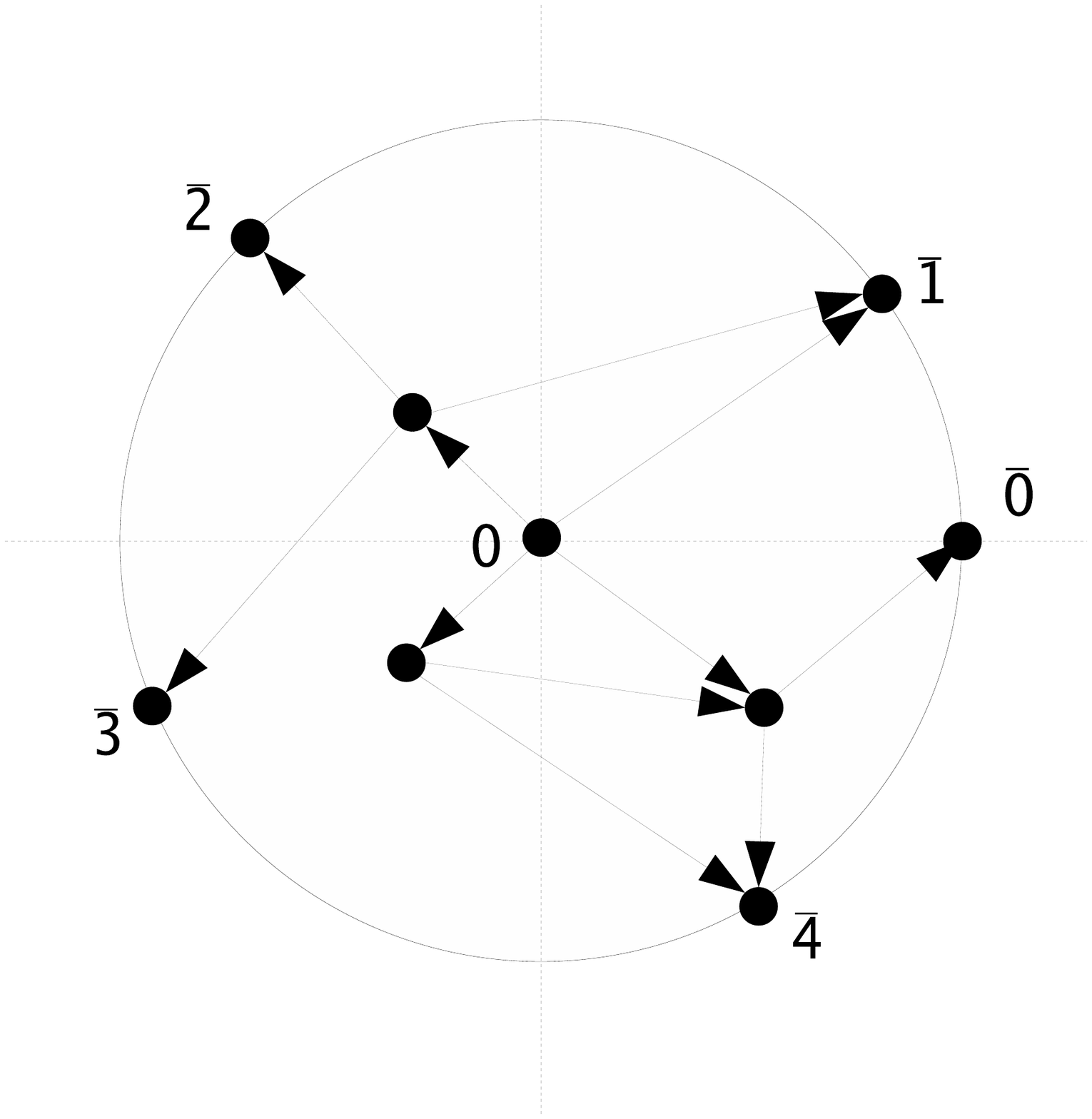}
\caption{\label{fig:exgraph} Some Shoikhet graph. The set of central edges $E_c(\Gamma)$ in this case consists of the four edges starting at the vertex ``$0$'' in the middle.}
\end{figure}

An example graph embedded in $D$ is shown in Figure \ref{fig:exgraph}.

The differential form $\omega_\Gamma$ that is integrated over configuration space can be expressed as a product of 
one-forms, one for each edge in $\Gamma$. 
\begin{align*}
\omega_\Gamma 
&=
\bigwedge_{(0\edge K)\in E_{c}(\Gamma)} d\theta_c(z_K, z_{\bar{0}}) \wedge \bigwedge_{j=1}^n \bigwedge_{(j\edge L)\in E(\Gamma)} d\theta(z_j,z_L)
\end{align*}

Here the one-forms occuring are defined as
\begin{align}
\label{equ:thetacdef}
d\theta_c(z, w) &= -\frac{1}{2 \pi}d\arg \left( \frac{z}{w} \right)\\ 
\label{equ:thetadef}
d\theta(z,w) &= \frac{1}{2 \pi} d\arg \left( (z-w)(1-z\bar{w})\bar{z} \right)
\end{align}

\begin{figure}
\includegraphics[trim = 0mm 40mm 0mm 20mm, clip, width=.4\textwidth]{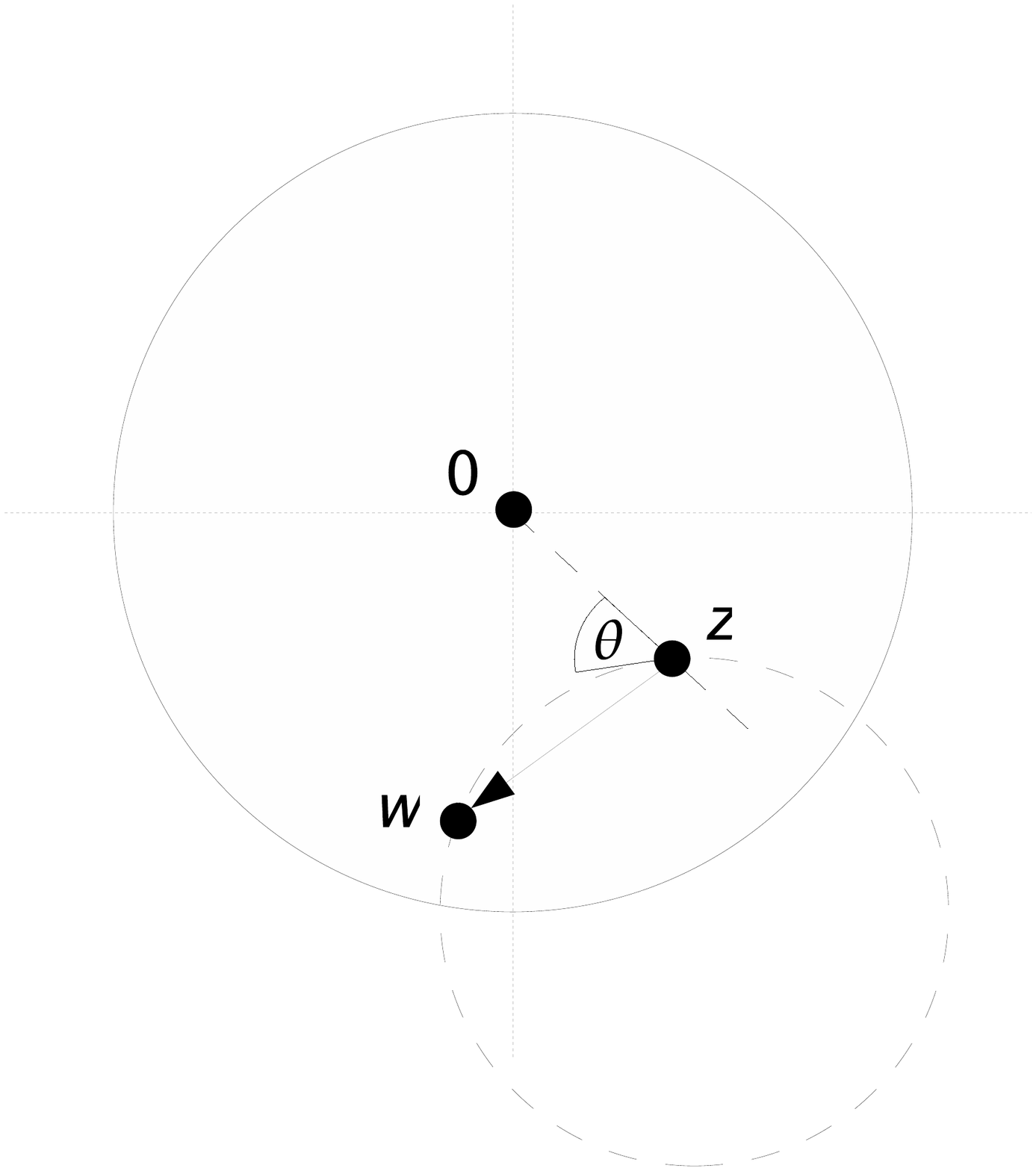}
\includegraphics[trim = 0mm 40mm 0mm 20mm, clip, width=.4\textwidth]{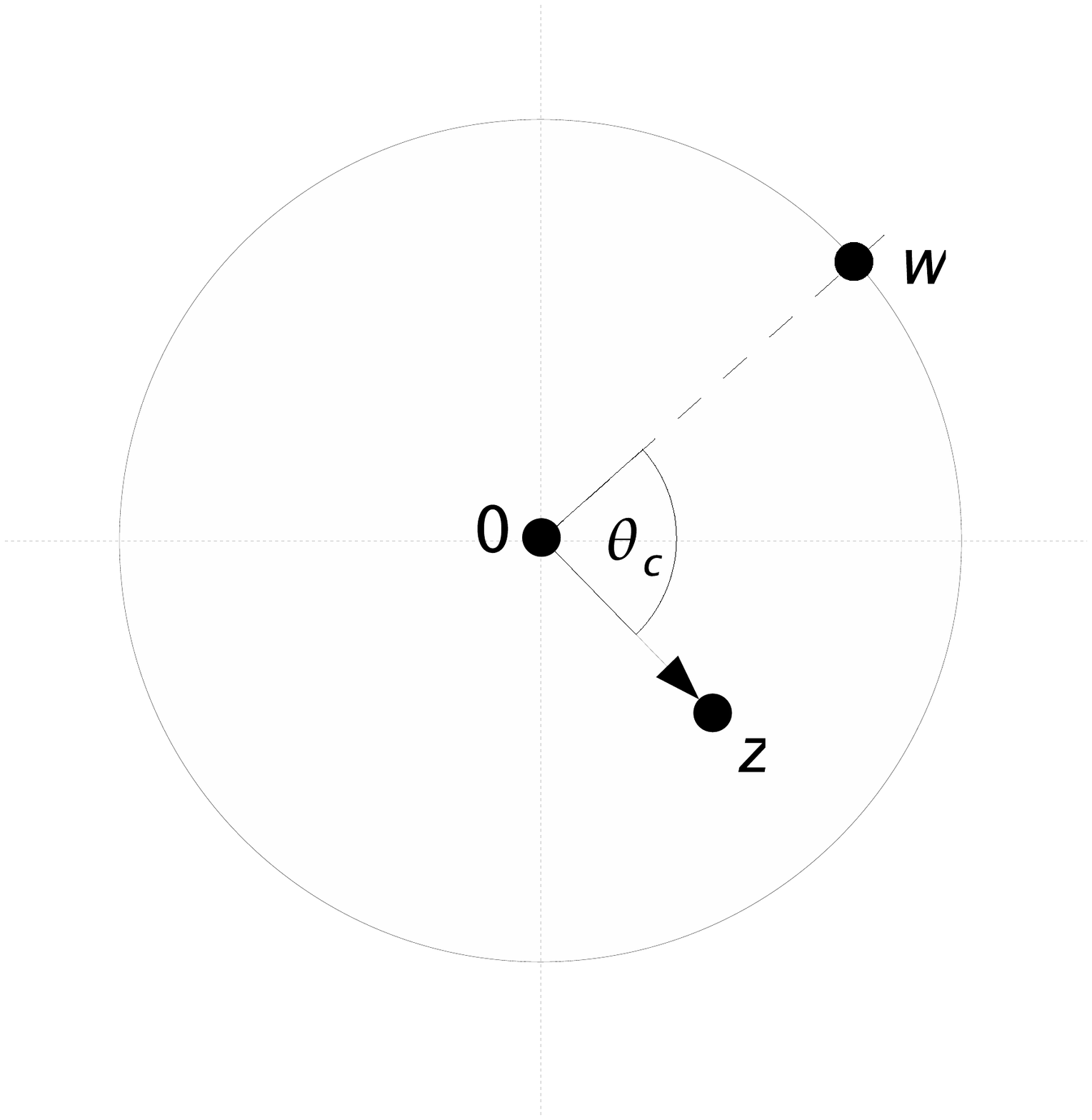}
\caption{\label{fig:phis} Geometric meaning of Shoikhet's angle forms.}
\end{figure}

The geometric meaning of these forms is illustrated in Figure \ref{fig:phis}. The ordering of the forms within the wedge products is such that forms corresponding to edges with source vertex $j$ stand on the left of those with source vertex $j+1$, and according to the order given on the stars for edges having the same source vertex.

We will use the abbreviations
\begin{align*}
\omega_\Gamma^c = \bigwedge_{(0\edge K)\in E_{c}(\Gamma)} d\theta_c(z_K, z_{\bar{0}})
& \quad \quad 
\omega_\Gamma^{nc} = \bigwedge_{j=1}^n \bigwedge_{(j\edge L)\in E(\Gamma)} d\theta(z_j,z_L)
\end{align*}
for the factors of $\omega_\Gamma = \omega_\Gamma^c \wedge \omega_\Gamma^{nc}$ coming from central and non-central edges.

\begin{rem}
All the differential forms above are defined on the enlarged configuration space $\tilde{C}_\Gamma$. The integral in the definition of the weights \eref{equ:wgammadef} shall be understood as the integral along the compact submanifold $C_\Gamma \subset \tilde{C}_\Gamma$ of the form $\omega_\Gamma$ on $\tilde{C}_\Gamma$.
\end{rem}

\begin{rem}
\label{rem:thetac}
Note that the form $d\theta_c(z, w)$ satisfies $d\theta_c(z, w) = d\theta_c(z, u) + d\theta_c(u, w)$ for any $u\in D\setminus\{0\}$.
\end{rem}

\subsection{Several remarks on orientations, signs, and rotation invariance}
On $\tilde{C}_\Gamma$ there is an obvious $S^1$-action by rotations, and $C_\Gamma$ intersects each $S^1$-orbit exactly once. Furthermore, note that the form $\omega_\Gamma$ is $S^1$-basic. In the following, fix a generator $\zeta$ of this action, generating a counterclockwise rotation. This is equivalent to choosing an orientation on $S^1$.

On $C_\Gamma$ we will then put the orientation that is induced by $\zeta$ and the volume form $\Omega = dz_1d\bar{z}_1\dots dz_md\bar{z}_m d\arg z_{\bar{0}}\dots d\arg z_{\bar{n}}$ on $\tilde{C}_\Gamma$.\footnote{This means, that the orientation on $C_\Gamma$ is determined by the form $\left.\iota_\zeta \Omega\right|_{C_\Gamma} = dz_1d\bar{z}_1\dots dz_md\bar{z}_md\arg z_{\bar{1}}\dots d\arg z_{\bar{n}}$. } 

Next consider the space $C_\Gamma^{(j)} = \{z_{\bar{j}}=1\} \subset \tilde{C}_\Gamma$, with the orientation determined by $\iota_\zeta \Omega$. It is not hard to show, using the homotopy by rotations of $C_\Gamma$ and $C_\Gamma^j$ and the fact that $\omega_\Gamma$ is $S^1$-basic\footnote{$S^1$-invariance would not be enough}, that
\begin{equation}
\label{equ:intsequal}
\int_{C_\Gamma }\omega_\Gamma=\int_{C_\Gamma^{(j)} }\omega_\Gamma.
\end{equation}

\section{Proof of Theorem \ref{thm:main}}
\label{sec:theproof}
We have to show that
\begin{equation}
\label{equ:prooftbs}
\left(d\morphU^{sh}_m(\gamma_1\wedge\dots \wedge\gamma_m)(a_0\otimes\dots\otimes a_n)\right)[\xi]
=
\morphU^{sh}_m(\gamma_1\wedge\dots \wedge\gamma_m)(B(a_0\otimes\dots\otimes a_n))[\xi].
\end{equation}
In fact, we will show that both sides of the above equation equal the following expression.
\begin{equation}
\label{equ:dashed}
\sum_{\Gamma\in G(n,m)} w_{\Gamma-\{e\}} D_\Gamma(\xi,\gamma_1\wedge\dots \wedge\gamma_m; a_0,..,a_n)
\end{equation}
Here $e$ ist the first edge in $E_c(\Gamma)=Star(0)$. 

\begin{lemma}
The l.h.s. of \eref{equ:prooftbs} is equal to \eref{equ:dashed}.
\end{lemma}
\begin{proof}
We can assume w.l.o.g. that $\xi=\xi_1\wedge\dots\wedge\xi_p$, with the $\xi_j$ constant vector fields. Then, for any form $\omega$, we have
\[
( d\omega )[\xi]= \sum_{i=1}^p (-1)^{i+1} \xi_i \cdot \omega[\xi_1\wedge ..\wedge \hat{\xi}_i\wedge ..\wedge \xi_p].
\]
On the other hand we have 
\begin{multline*}
\sum_{i=1}^p (-1)^{i+1} \xi_i \cdot
D_\Gamma(\xi_1\wedge ..\wedge \hat{\xi}_i\wedge ..\wedge \xi_p,\gamma_1\wedge\dots \wedge\gamma_m; a_0,..,a_n)
= \\
\sum_{v\in V(\Gamma) \setminus \{0\}} D_{\Gamma\cup \{(0\edge v)\}}(\xi, \gamma_1\wedge\dots \wedge\gamma_m, a_0,..,a_n)
\end{multline*}
Here by $\Gamma\cup \{(0\edge v)\}$ we mean the graph formed by adding the edge $(0\edge v)$ to $\Gamma$ and adjusting the ordering in $E_c(\Gamma)$ so that the newly added edge is the first.
Next multiply by $w_\Gamma$ and sum over all graphs $\Gamma$. Observe that the double sum occuring, namely
\[
\sum_{\Gamma\in G(n,m)} \sum_{v\in V(\Gamma) \setminus \{0\}}
\]
contains every graph in $G(n,m)$ (i.e., a graph with ordering on the stars) exactly once. Hence the Lemma has been shown.
\end{proof}

\begin{lemma}
\label{lem:rhsequal}
The r.h.s. of \eref{equ:prooftbs} is equal to \eref{equ:dashed}.
\end{lemma}
For the proof, we need some preparation. First define the operator $\sigma$ (cyclic shift) on $C_\bullet(A,A)$ by
\[
\sigma(a_0\otimes a_1\otimes \dots \otimes a_n)= a_0\otimes a_2\otimes a_3\otimes\dots\otimes a_n\otimes a_1.
\]
Similarly define the cyclic shift operator, also called $\sigma$, on a graph $\Gamma$ by cyclically interchanging the labels on the type II vertices except the vertex $\bar{0}$, such that the following holds
\[
D_{\sigma\Gamma}(\dots,a_0\otimes \dots \otimes a_n)=D_{\Gamma}(\dots,\sigma(a_0\otimes \dots \otimes a_n)).
\]
Also define the operator $s$ on $C_\bullet(A,A)$ by
\[
s(a_0\otimes a_1\otimes \dots \otimes a_n) = 1\otimes a_0\otimes a_1\otimes \dots \otimes a_n
\]
so that $B=\sum_{i=0}^n (-1)^{in} \sigma^i s$. We can also define the operator $s$ on graphs so that
\[
D_{s\Gamma}(\dots,a_0\otimes \dots \otimes a_n)=D_{\Gamma}(\dots,s(a_0\otimes \dots \otimes a_n)).
\]
The graph $s\Gamma$ is the same as $\Gamma$ but with the vertex $\bar{0}$ deleted and the type II vertices renumbered such that $\bar{1}\in V(\Gamma)$ becomes $\bar{0}\in V(s\Gamma)$, $\bar{2}\in V(\Gamma)$ becomes $\bar{1}\in V(s\Gamma)$ etc. In case there is an edge in $\Gamma$ ending at $\bar{0}$, we will set $s\Gamma = \emptyset$ the empty graph and define $D_{s\Gamma}(\dots):=0$.

Then we can compute
\begin{align}
\sum_\Gamma w_\Gamma D_{\Gamma}(\dots,B(a_0\otimes \dots \otimes a_n)) 
&=
\sum_{i=0}^{n} (-1)^{in} \sum_\Gamma w_\Gamma D_{\Gamma}(\dots,\sigma^i(1\otimes a_0\otimes \dots \otimes a_n)) \nonumber \\
&=
\sum_{i=0}^{n} (-1)^{in} \sum_\Gamma w_\Gamma D_{\sigma^i\Gamma}(\dots,1\otimes a_0\otimes \dots \otimes a_n) 
\label{equ:prerhsequal} \\
&=
\sum_{i=0}^{n} (-1)^{in} \sum_\Gamma w_{\sigma^{-i}\Gamma} D_{\Gamma}(\dots,1\otimes a_0\otimes \dots \otimes a_n) \nonumber \\
&=
\sum_\Gamma \left( \sum_{i=0}^{n} (-1)^{in} w_{\sigma^{-i}\Gamma}\right)  D_{\Gamma}(\dots,1\otimes a_0\otimes \dots \otimes a_n) \nonumber
\end{align}
where in the second to last equality we changed variables in the $\Gamma$-summation.\footnote{We used that $\sigma$ is a bijection on the set of graphs $G(n,m)$.}

\begin{lemma}
\label{lem:weightsid}
Let $\Gamma$ be a graph with $n+2$ type II vertices and no edge hitting the $\bar{0}$ vertex. Then
\[
\sum_{i=0}^{n} (-1)^{in} w_{\sigma^{-i}\Gamma} = \frac{1}{\#E_c(s\Gamma)} \sum_{i=1}^{\# E_c(s\Gamma)} (-1)^{i+1} w_{s\Gamma-\{e_i\}}.
\]
where $e_i$ is the $i$-th edge in $E_c(s\Gamma)$.
\end{lemma}
\begin{proof}
We will show the equality from right to left. 

\begin{align}
&\frac{1}{\#E_c(s\Gamma)}
\sum_{i=1}^{\# E_c(s\Gamma)} (-1)^{i+1} w_{s\Gamma-\{e_i\}}= \nonumber \\
&\quad \quad=
\prod_{v\in V(\Gamma)} \frac{1}{(\#Star(v))!}
\sum_{i=1}^{\# E_c(s\Gamma)} (-1)^{i+1}
\int_{C_{s\Gamma}} 
\bigwedge_{(0\edge K)\in E_{c}(s\Gamma-\{e_i\})} d\theta_c(z_K,z_{\bar{0}}) \wedge \omega_{s\Gamma}^{nc}  \nonumber \\
&\quad \quad=
\prod_{v\in V(\Gamma)} \frac{1}{(\#Star(v))!}
\int_{C_{s\Gamma}} \int_{Z\in S^1}
\bigwedge_{(0\edge K)\in E_{c}(s\Gamma)} (d\theta_c(z_K,z_{\bar{0}})+d\theta_c(z_{\bar{0}},Z)) \wedge \omega_{s\Gamma}^{nc} \label{equ:lemprerhs1} \\
&\quad \quad=
\prod_{v\in V(\Gamma)} \frac{1}{(\#Star(v))!}
\int_{C_{s\Gamma}} 
\sum_{i=0}^n \int_{\arg{Z}\in (\arg z_{\widebar{i-1}},z_{\widebar{i}})}
\bigwedge_{(0\edge K)\in E_{c}(s\Gamma)} d\theta_c(z_K,Z) \wedge \omega_{s\Gamma}^{nc}.
\nonumber
\end{align}

In the last line we used Remark \ref{rem:thetac} and (independently) decomposed the domain of the $Z$-integral into $n+1$ pieces. Consider only the $i-th$ piece:
\[
\int_{U_i\cap \{z_{\bar{0}}=1\}} \bigwedge_{(0\edge K)\in E_{c}(s\Gamma)} d\theta_c(z_K,Z) \wedge \omega_{s\Gamma}^{nc}
\]
where $U_i\subset S^1 \times \tilde{C}_{s\Gamma}$ is the set 
\[
U_i = \{(Z, z_0,..,z_m,z_{\bar{0}},..,z_{\bar{n}})\in S^1\times \tilde{C}_{s\Gamma} \mid \arg{Z}\in (\arg z_{\widebar{i-1}},z_{\widebar{i}}) \}.
\]
Note that the set $U_i$ can be identified with an open dense subset of the enlarged configuration space $\tilde{C}_{\sigma^{-i}\Gamma}$ of the graph $\sigma^{-i}\Gamma$. Concretely, the map is given by
\begin{align*}
\phi_i : U_i &\to \tilde{C}_{\sigma^{-i}\Gamma} \\
(Z, z_0,..,z_m,z_{\bar{0}},..,z_{\bar{n}}) &\mapsto (z_0,..,z_m,Z, z_{\widebar{i+1}}, \dots ,z_{\bar{i}}).
\end{align*}

This map is in general not orientation preserving, but changes the orientation by a factor $(-1)^{in}$. Note also that 
\[
\phi_i^* \omega_{\sigma^{-i}\Gamma}=\bigwedge_{(0\edge K)\in E_{c}(s\Gamma)} d\theta_c(z_K,Z) \wedge \omega_{s\Gamma}^{nc}
\]
is exactly the differential form integrated over in the above integral. Hence the $i$-th piece considered above can be rewritten as 
\begin{align*}
&\int_{U_i\cap \{z_{\bar{0}}=1\}} \bigwedge_{(0\edge K)\in E_{c}(s\Gamma)} d\theta_c(z_K,Z) \wedge \omega_{s\Gamma}^{nc}  =
\int_{U_i\cap \{z_{\bar{0}}=1\}} \phi_i^* \omega_{\sigma^{-i}\Gamma} \\
&\quad=
(-1)^{in}
\int_{\{z_{\widebar{I(i)}}=1\} \subset \tilde{C}_{\sigma^{-i}\Gamma} } \omega_{\sigma^{-i}\Gamma}
=
(-1)^{in}
\int_{C_{\sigma^{-i}\Gamma}} \phi_i^* \omega_{\sigma^{-i}\Gamma}.
\end{align*}

In the second line, the function $I(i)$ is defined such that $I(0)=1$, $I(1)=n+1$, $I(2)=n$ etc. Put differently, it is
defined such that the coordinate function $z_{\bar{0}}$ on $U_i$ is the pullback $\phi_i^* z_{\widebar{I(i)}}$ of the coordinate function $z_{\widebar{I(i)}}$ on $\tilde{C}_{\sigma^{-i}\Gamma}$.\footnote{The author apologizes for using the same symbol $z_{\bar{j}}$ for two different functions on two different spaces. However, adding a superscript indicating the space would make the notation rather clumsy.} In the last line we furthermore used \eref{equ:intsequal}.

Inserting into \eref{equ:lemprerhs1} we finally obtain
\begin{align*}
\frac{1}{\#E_c(s\Gamma)}
\sum_{i=1}^{\# E_c(s\Gamma)} (-1)^{i+1} w_{s\Gamma-\{e_i\}} 
&=
\prod_{v\in V(\Gamma)} \frac{1}{(\#Star(v))!}
\sum_{i=0}^n (-1)^{in}\int_{C_{\sigma^{-i}\Gamma}} \omega_{\sigma^{-i}\Gamma} \\
&= 
\sum_{i=0}^{n} (-1)^{in} w_{\sigma^{-i}\Gamma}
\end{align*}

\end{proof}

With this Lemma at hand, we can now finish the proof of Lemma \ref{lem:rhsequal}.

\begin{proof}[Proof of Lemma \ref{lem:rhsequal}]
Continue the computation \eref{equ:prerhsequal}. We get, using the previous Lemma
\begin{align*}
& \sum_{i=0}^{n} (-1)^{in} \sum_{\Gamma \in G(n,m+1)} w_{\sigma^{-i}\Gamma} D_{\Gamma}(\dots,1\otimes a_0\otimes \dots \otimes a_n) \\
&=
\sum_{\Gamma \in G(m,n+1)} \frac{1}{\#E_c(s\Gamma)}\sum_{i=1}^{\# E_c(s\Gamma)} (-1)^{i+1}w_{s\Gamma-\{e_i\}} D_{\Gamma}(\dots,1\otimes a_0\otimes \dots \otimes a_n) \\
&=
\sum_{\Gamma \in G(m,n+1)} \frac{1}{\#E_c(s\Gamma)}\sum_{i=1}^{\# E_c(s\Gamma)} (-1)^{i+1}w_{s\Gamma-\{e_i\}} D_{s\Gamma}(\dots,a_0\otimes \dots \otimes a_n) \\
&=
\sum_{\Gamma \in G(m,n)} \frac{1}{\#E_c(\Gamma)}\sum_{i=1}^{\# E_c(\Gamma)} (-1)^{i+1}w_{\Gamma-\{e_i\}} D_{\Gamma}(\dots,a_0\otimes \dots \otimes a_n)  \\
&= 
\sum_{\Gamma' \in G(m,n)} w_{\Gamma'-\{e\}} D_{\Gamma'}(\dots,a_0\otimes \dots \otimes a_n)  \\
\end{align*}
For the second to last equality, we used that the map $s$ is surjective and changed variables. In the last line $e$ is again the first edge of $E_c(\Gamma')$. For the last equality we also ``changed variables''. We replaced each pair $(\Gamma,i)$ by the pair $(\Gamma', i)$, where $\Gamma'$ is the same graph as $\Gamma$, but with the ordering on $Star(0)$ changed by putting the $i$-th edge at first position in the ordering. Then 
\[
w_{\Gamma-\{e_i\}} = w_{\Gamma'-\{e\}}
\]
and 
\[
(-1)^{i+1} D_\Gamma(\dots) =  D_{\Gamma'}(\dots).
\]
In the resulting sum, everything is independent of $i$, and the $i$-summation just cancels the factor $\frac{1}{\#E_c(\Gamma)}$. Hence the lemma and thus Theorem \ref{thm:main} has been proven.
\end{proof}

\appendix

\section{Standard Definitions, Gradings and Signs}
In this section, we recite some standard definitions and results. We mostly use the terminology of Tsygan \cite{tsygan}, and hence almost copy the expositions given in his paper.
\subsection{$L_\infty$-algebras and $L_\infty$-modules}
Let $\alg{g}^\bullet$ be a $\mathbb{Z}$-graded vector space. An $L_\infty$-structure on $\alg{g}^\bullet$ is a degree $1$ coderivation $Q$ on the cocommutative cofree coalgebra $S(\alg{g}^\bullet [1])$ satisfying 
\[
Q^2=0.
\]
Any coderivation on $S(\alg{g}^\bullet [1])$ is determined by its projection to $\alg{g}^\bullet$, hence by a series of linear functions
\[
q_k \in Hom(\bigwedge^k \alg{g}^\bullet, \alg{g}^\bullet)
\]
of degree $2-k$. The condition that $Q^2=0$ reads
\[
\sum_{j=1}^N \sum_{\sigma \in S_N}  \pm \frac{1}{j!(N-j)!} q_{N-j+1}(q_{j}(a_{\sigma(1)},..,a_{\sigma(j)}),a_{\sigma(j+1)}),..,a_{\sigma(N)})) = 0
\]
for all $N=1,2,..$ and all $a_1,..,a_N\in \alg{g}^\bullet$. Here the sign is the lexicographic sign w.r.t. the shifted-by-one grading.

Let now $M^\bullet$ be another graded vector space. An $L_\infty$-module structure on $M^\bullet$ is a degree $1$ coderivation $D$ on the free comodule 
\[
S(\alg{g}^\bullet [1])\otimes M^\bullet
\]
satisfying $D^2=0$.
Again, $D$ is determined by its composition with the projection to $M^\bullet$, i.e., by components
\[
d_k \in Hom(\bigwedge^k \alg{g}^\bullet\otimes M^\bullet, M^\bullet)
\]
of degree $1-k$ such that the following holds for all $N=1,2,..$ and $a_1,..,a_N\in \alg{g}^\bullet, m\in M^\bullet$:
\begin{multline*}
\sum_{j=1}^N \sum_{\sigma \in S_N} \left[ 
\pm  \frac{1}{j!(N-j)!}
d_{N-j}(a_{\sigma(1)},..,a_{\sigma(j)},d_j(a_{\sigma(j+1)},..,a_{\sigma(N)},m)) \right.
\\
\left.
\pm
\frac{1}{j!(N-j)!} d_{N-j+1}(q_{j}(a_{\sigma(1)},..,a_{\sigma(j)}),a_{\sigma(j+1)},..,a_{\sigma(N)},m) \right] = 0
\end{multline*}

Morphisms of $L_\infty$-algebras and $L_\infty$-modules are defined in the obvious way as morphisms of the underlying coalgebras or comodules that commute with the structure ($Q$ or $D$) given.

Philosophically, and also mathematically if $\dim \alg{g}^\bullet<\infty$, one can understand the components $q_k$ of $Q$ as terms in a ``Taylor series'' 
\[
Q = \sum_{k\geq 1} \frac{q_k}{k!}
\]
of a degree $1$ vector field $Q$ on $\alg{g}^\bullet[1]$, commuting with itself. Consider next the trivial bundle $\alg{g}^\bullet [1]\otimes M^\bullet \rightarrow \alg{g}^\bullet [1]$. An $L_\infty$-module structure can be understood philosophically as a flat lift $D$ of the vector field $Q$ to this bundle.

\begin{rem}
\label{rem:whattoproof}
The only way in which the above definitions are needed in this paper is the following. Consider an $L_\infty$-algebra $(\alg{g}^\bullet,Q)$ as above and a morphism $\morphU$ of $L_\infty$-modules over $\alg{g}^\bullet$
\[
\morphU: (M_1^\bullet,D_1) \rightarrow (M_2^\bullet,D_2).
\]
We next want to modify the $L_\infty$-module structures to 
\begin{align*}
D_1' &= D_1 + \delta_1 \\
D_2' &= D_2 + \delta_2 
\end{align*}
where the $\delta_j$ are degree $1$ endomorphisms of $S\alg{g}^\bullet[1]\otimes M_j^\bullet$.
Then $\morphU$ is still a morphism of the new $L_\infty$-modules $(M_j^\bullet,D_j')$ if and only if 
\[
\morphU \circ \delta_1 = \delta_2 \circ \morphU.
\]
As usual, it is sufficient to consider the projection of both sides to $M_2^\bullet$, because $D_j'$ are coderivations. In our case furthermore, all Taylor components of the $\delta_j$ vanish except in degree $0$. Hence the above condition reads in components
\[
\morphU_N(a_1,..,a_N,\delta_1 m)= \delta_2 \morphU_N(a_1,..,a_N,m)
\]
for $N=0,1,..$. This is precisely the condition \eref{equ:prooftbs} proven in Section \ref{sec:theproof}.
\end{rem}

\subsection{Polyvector Fields}
The grading we use on the space of polyvector fields $T_{poly}^\bullet$ is such that a vector field has degree $0$, a bivector field degree $1$, a function degree $-1$ etc. The Schouten-Nijenhuis bracket $\co{\cdot}{\cdot}_{SN}$ on $T_{poly}^\bullet$ is defined such that 
\begin{align*}
\co{f}{g}_{SN} &= 0 \\
\co{\xi}{\gamma_1}_{SN} &= L_\xi \gamma_1 \\
\co{\gamma_1}{\gamma_2\wedge \gamma_3}_{SN} &= \co{\gamma_1}{\gamma_2}_{SN} \wedge \gamma_3 +(-1)^{|\gamma_1|(|\gamma_2|+1)} \gamma_2\wedge \co{\gamma_1}{\gamma_3}_{SN}
\end{align*}

for all functions $f\in A$, vector fields $\xi\in T_{poly}^0$ and polyvector fields $\gamma_1,\gamma_2,\gamma_3\in T_{poly}^\bullet$. Note that the sign is the lexicographic one if we count $\wedge$ to have degree $+1$. This is as expected since 
\[
\wedge : T_{poly}^\bullet\otimes T_{poly}^* \rightarrow T_{poly}^{\bullet+*+1}.
\]

One can check that the above bracket turns $T_{poly}^\bullet$ into a graded Lie algebra. As any Lie algebra, it is automatically an $L_\infty$-algebra, obtained by setting
\[
q_k = 
\begin{cases}
\co{\cdot}{\cdot}_{SN} & \quad \text{for $k=2$} \\
0 & \quad \text{otherwise}.
\end{cases}
\]

Next consider the space $\Omega^\bullet(M)$ of differential forms on the manifold $M$. We consider it with the opposite of the usual grading, i.e., a $k$-form has degree $-k$. With this grading, $\Omega^\bullet(M)$ is a graded module over the graded Lie algebra $T_{poly}^\bullet$. The action is given by 
\[
\gamma\otimes \omega \rightarrow L_\gamma \omega = \co{d}{\iota_\gamma} \omega
\]
for polyvector fields $\gamma$ and differential forms $\omega$. For a function $f\in T_{poly}^1$ we define $\iota_f$ to be the multplication by $f$. Any module over a Lie algebra is also an $L_\infty$-module, in this case by setting 
\[
d_k(\gamma_1,..,\gamma_k,\omega) =
\begin{cases}
L_{\gamma_1}\omega & \quad \text{for $k=1$} \\
0 & \quad \text{otherwise}.
\end{cases}
\]

\subsection{Hochschild and Cyclic Cohomology}
The Hochschild cochain complex $C^\bullet(A,A)$ of the unital algebra $A$ with values in the $A$-bimodule $M$ is defined as
\[
C^k(A,M) = \Hom(A^{\otimes k},M).
\]
The Hochschild coboundary operator $d_H$ is given by
\begin{multline*}
(d_H\Psi)(a_1,..,a_{n+1}) = (-1)^{n+1} a_1\Psi(a_2,..,a_{n}) +\\
+ \sum_{j=1}^n (-1)^{j+n+1} \Psi(a_1,..,a_{j-1},a_j a_{j+1},a_{j+2},..,a_{n+1} )
+ \Psi(a_1,..,a_n) a_{n+1}.
\end{multline*}

There is a Lie bracket $\co{\cdot}{\cdot}_G$ on $C^\bullet(A,A)[1]$, called the Gerstenhaber bracket. it is defined as \[
\co{\Psi}{\Phi}_G = \Psi \circ \Phi - (-1)^{(m-1)(n-1)} \Phi\circ \Psi
\]
where $\Psi\in C^m(A,A),\Phi\in C^n(A,A)$ and 
\begin{multline*}
(\Psi\circ \Phi) (a_1,..,a_{n+m-1}) \\
= \sum_{j=1}^{m} (-1)^{(n-1)(j-1)} \Psi(a_1,..,a_{j-1},\Phi(a_j,..,a_{j+n-1}),a_{j+m},..,a_{n+m-1}).
\end{multline*}

If we set 
\[
m(a1,a_2) = a_1\cdot a_2,
\]
so that $m\in C^1(A,A)$, one can check that $d_H(\cdot) = \co{m}{\cdot}_G$. Hence, by the Jacobi identity for $\co{\cdot}{\cdot}_G$, $C^\bullet(A,A)[1]$ is a differential graded Lie algebra, and hence an $L_\infty$-algebra. 

The normalized Hochschild chain complex $C_\bullet(A,M)$ with values in the bimodule $M$ is defined as
\[
C_k(A,M) = M\otimes \bar{A}^{\otimes k}
\]
where $\bar{A}= A/(1\cdot \cn{})$. The differential is
\begin{multline*}
b(m\otimes a_1\otimes \cdots \otimes a_n) = m\cdot a_1\otimes a_2\otimes\cdots \otimes a_n + \\
+\sum_{j=1}^{n-1}(-1)^j m\otimes a_1\otimes\cdots \otimes a_ja_{j+1}\otimes \cdots \otimes a_n
+ (-1)^n a_n\cdot m\otimes a_1\otimes\cdots\otimes a_{n-1}.
\end{multline*}

The action \eref{equ:hochaction} makes $C_\bullet(A,A)$ with the opposite (negative) grading into a differential graded module over $C^\bullet(A,A)[1]$. On $C_\bullet(A,A)$ there is another natural operation, namely the $B$ of \eref{equ:Bdef}. One can check that $B$ anticommutes with $b$, so that it makes sense to define the cyclic chain complex $(CC^W_\bullet(A,A), b+uB)$ as in \eref{equ:cyccdef}. Depending on the choice of the $\cn{}((u))$-module $W$ one obtains different cyclic cohomology theories:

\begin{itemize}
\item For $W=\cn{}$ with $u$ acting as $0$ one recovers the usual Hochschild chain complex.
\item For $W=\C((u))$ one obtains the periodic cyclic chain complex $CC^{per}_\bullet(A,A)$. In the case $A=C^\infty(M)$, it is isomorphic to the complex $(\Omega^\bullet(M) ((u)),d)$, whose cohomology is $H^\bullet(M) ((u))$.
\end{itemize}

Furthermore $B$ (graded) commutes with the action of $C^\bullet(A,A)[1]$, and hence the cyclic chain complex carries the structure of a differential graded $C^\bullet(A,A)[1]$-module.

In the case of interest to us, the algebra $A=C^\infty(M)$ is a locally convex algebra, and the tensor products occuring in the above definitions shall be understood as projectively completed tensor products (see \cite{connes}, section 5). 

\nocite{*}
\bibliographystyle{plain}
\bibliography{tsygancyc} 

\end{document}